\documentclass[12pt]{article}

\usepackage{graphicx,tikz,color,url}
\usepackage{amsmath,amssymb}
\usepackage{cite}
\usetikzlibrary{arrows}
\usepackage{float}
\usepackage{tikz}
\usepackage{arydshln}
\usepackage{tikz}
\usetikzlibrary{calc}
\usepackage[labelformat=simple]{subcaption}

\usetikzlibrary{positioning}
\usepackage{changepage}
\usepackage{pict2e}
\usepackage{chemfig}

\newtheorem{theorem}{Theorem}
\newtheorem{definition}{Definition}
\newtheorem{lemma}{Lemma}
\begin{document}

\allowdisplaybreaks

\title{Enumerating the number of $k$-matchings in successively amalgamated graphs}

\author{Simon Grad$^{a,}$\thanks{email: \texttt{sg12625@student.uni-lj.si}} 
\and Sandi Klav\v zar$^{a,b,c,}$\thanks{corresponding author; email: \texttt{sandi.klavzar@fmf.uni-lj.si}}
}
\maketitle

\begin{center}
$^a$ Faculty of Mathematics and Physics, University of Ljubljana, Slovenia\\
\medskip

$^b$ Institute of Mathematics, Physics and Mechanics, Ljubljana, Slovenia\\
\medskip
	
$^c$ Faculty of Natural Sciences and Mathematics, University of Maribor, Slovenia\\
\medskip

\end{center}
\begin{abstract}
In this paper, the transfer matrix technique using the $k$-matching vector is developed to compute the number of $k$-matchings in an arbitrary graph  which can be constructed by successive amalgamations over sets of cardinality two. This widely extends known methods from the literature developed for computing the number of $k$-matchings in benzenoid chains, octagonal chains, cyclooctatetraene chains, and arbitrary cyclic chains. Two examples demonstrating how the present method can be applied are given, one of then being an elaborated chemical example.
\end{abstract}


\noindent
{\bf AMS Subj.\ Class.\ (2020)}: 05C92, 05C09, 05C70 

%%%%%%%%%%%%%%%%%%%%%%%%%%%%%%
\section{Introduction}
%%%%%%%%%%%%%%%%%%%%%%%%%%%%%%

Matchings form a central graph concept that is also relevant and useful in many other fields. From our point of view, the most relevant application of matchings is in chemistry, since matchings naturally manifest double bonds in chemical compounds. In chemistry, matchings are also known as Kekul\'e{} structures, see the book~\cite{cyvin-1988} solely devoted to the Kekul\'e{} structures in benzenoid hydrocarbons. As a rule of thumb, the more matchings a molecular graph has, the more stable the molecule is. 

From the above reasons, the enumeration of $k$-matchings is a fundamental problem. In particular, the number of perfect matchings of a molecular graph quantifies aromaticity and aids in analyzing resonance energy and chemical stability~\cite{joice-2023}. The total number of matchings in a graph is called the Hosoya index of a graph, sometimes also called Z-index. It was introduced in~\cite{hosoya-1971}, and studied from then on in different directions. We refer to the 2010 survey article~\cite{wagner-2010} and the following selected subsequent papers~\cite{cruz-2022, doslic-2024, huang-2018, kuzmin-2024, oo-2024, sahin-2022, vesel-2021}. 

In a couple of papers~\cite{polansky-1989, randic-1989}, Polansky, Randi\'c, and Hosoya, used what is now known as the transfer matrix technique to compute the Hosoya index, the matching polynomial, and the Wiener index of benzenoid chains. The idea is to assign a vector to a benzenoid system such that with successive multiplications with a transfer matrix (reflecting additions of hexagons) the value of the desired invariant is obtained. 

To compute the Hosoya index of catacondensed hexagonal systems, Cruz et al.~\cite{cruz-2017} invented the so-called Hosoya vector and adopted the transfer matrix method for this situation. Oz and Cangul~\cite{oz-2022} followed by introducing the so-called $k$-matching vector in order to get the number of $k$-matchings in benzenoid chains. The approach has been extended to  cyclooctatetraene chains~\cite{chen-2024}, to arbitrary cyclic chains~\cite{grad-2025}, to branched benzenoid systems~\cite{oz-2023}, and elsewhere~\cite{oz-2022b, farooq-2024, shi-2023}. In this paper, we develop the $k$-matching vector based transfer matrix method for arbitrary graphs which can be build up from smaller ingredients using amalgamations over vertex sets of cardinality two. In this way, we broadly generalize the approaches described above. 

In the following section we first introduce key concepts and review existing results. Then we define a family of graphs ${\cal U}$ which, roughly speaking, contains graphs that alow a decomposition into two subgraphs that are amalgamated over a vertex sets of cardinality two. We end the section by formulating our main result (Theorem~\ref{thm:main}) and discuss how it can be applied. Theorem~\ref{thm:main} is proved in Section~\ref{sec:proof}. In the final section we give two examples demonstrating how the developed method can be applied, the second one being an elaborated chemical example. 

%%%%%%%%%%%%%%%%%%%%%%%%%%%%%%
\section{Known and new concepts, the main theorem}
\label{sec:preliminaries}
%%%%%%%%%%%%%%%%%%%%%%%%%%%%%%

A {\em matching} of a graph $G=(V(G), E(G))$ is a set of edges $M\subseteq E(G)$ such that no two edges from $M$ are incident. We further say that $M$ is a {\em $k$-matching} if $|M| = k$. The total number of different matchings of $G$ is the {\em Hosoya index} of $G$, denoted by $Z(G)$. Setting $p(G,k)$ to denote the number of $k$-matchings of $G$, the Hosoya index can be written as follows: 
$$Z(G) = \sum_{k\ge 0} p(G,k)\,.$$   

In the following, we will use explicitly and implicitly the following well-known, simple lemmas. 

\begin{lemma}
  \label{lemma1}
If $ab$ is an edge of a graph $G$, then 
$$p(G,k) = p(G-ab,k) + p(G-a-b, k-1)\,.$$
\end{lemma}

\begin{lemma}
\label{lemma2}
If $G$ is the disjoint union of connected graphs $G_1$ and $G_2$, then 
$$p(G,k) = p(G_1 \cup G_2, k) = \sum_{i=0}^k p(G_1,i)p(G_2,k-i)\,.$$
\end{lemma}

The transfer matrix method is based on the following concept. If $a$ and $b$ are vertices in a graph $G$, then the {\em $k$-matching vector} $p_{ab}(G,k)$ of $G$ with respect to $ab$ is the vector 

$$p_{ab}(G,k) =
    \begin{pmatrix}
     p(G,k) \\
     p(G,k-1) \\
     \vdots \\
     p(G, 0) \\
     p(G - a, k) \\
     p(G - a, k-1) \\
     \vdots \\
     p(G - a, 0) \\
     p(G - b, k) \\
     p(G - b, k-1) \\
     \vdots \\
     p(G - b, 0) \\
     p(G - a - b, k) \\
     p(G - a - b, k-1) \\
     \vdots \\
     p(G - a - b, 0) \\
  \end{pmatrix}\,. $$

Based on the concept of the $k$-matching vector, for our purposes the following definition is crucial.

\begin{definition}[Family of graphs $\mathcal{U}$]
\label{def:key}
A graph $G$ belong to a family $\mathcal{U}$ if there exists a proper subgraph $G'$ of $G$, and vertices $x,y,a,b$, such that 
\begin{itemize}
\item $x,y\in V(G')$ and $G - \{x,y\}$ has precisely two components, and 
\item $\{a,b\}\ne \{x,y\}$. 
\end{itemize}
\end{definition}

The main result of this paper now reads as follows. 

\begin{theorem}
\label{thm:main}
If $G\in \mathcal{U}$, and $G'$, $x$, $y$, $a$, and $b$ are as in Definition~\ref{def:key}, then there exists a {\em transfer matrix} $A$ such that 
\begin{align*}
p_{ab}(G,k) & = A \cdot p_{xy}(G',k)\,.
\end{align*}
\end{theorem}

\begin{figure} [ht!]
    \centering
    \begin{tikzpicture}[scale=1.0,style=thick,x=1cm,y=1cm]
\def\vr{5pt}

% vertices defined
\coordinate(v1) at (0.0,1.5);
\coordinate(v2) at (0.0,-1.5);
\coordinate(v3) at (-2.0,1.0);
\coordinate(v4) at (-2.0,-1.0);
\coordinate(v5) at (2.0,1.0);
\coordinate(v6) at (2.0,-1.0);
\coordinate(v7) at (-4.0,1.0);
\coordinate(v8) at (-4.0,-1.0);
\coordinate(v9) at (3.0,2);
\coordinate(v10) at (4.0,1.0);
\coordinate(v11) at (4.0,-1.0);
\coordinate(v12) at (-5, 1.5);

%Inner duals
\coordinate(a1) at (0.0,0.0);
\coordinate(a2) at (3.0,0.0);
\coordinate(a3) at (-3.0,0.0);

% \edges                             
\draw (v3) -- (v4) -- (v8) -- (v7) -- (v3); 
\draw (v6) -- (v11) -- (v10) -- (v9) -- (v5);
\draw (v3) -- (v1) -- (v4);
\draw (v2) -- (v5); 
\draw (v6) -- (v2);
\draw (v9) -- (v6); 
\draw (v9) -- (v11); 
\draw (v5) -- (v10); 
\draw (v6) -- (v1); 
\draw (v5) -- (v10); 
\draw (v4) -- (v6); 
\draw (v4) -- (v2); 
\draw (v3) -- (v5); 
\draw (v1) -- (v5); 
\draw (v1) -- (v4); 
\draw (v2) -- (v3); 
\draw (v12) -- (v7);

\draw[rounded corners, densely dotted](0,1.8)--(-2.3,1.3)--(-2.3,-1.3)--(0,-1.8)--(2.3,-1.3)--(2.3,1.3)--cycle;
\draw[rounded corners, densely dotted](-1.7,-1.4)--(-1.7,1.4)--(-5.3,1.9)--(-4.3,-1.4)--cycle;
\draw[rounded corners, densely dotted](1.7,-1.4)--(1.7,2.3)--(4.3,2.3)--(4.3,-1.4)--cycle;

%  vertices
\foreach \i in {1,2,3,4,5,6,7,8,9,10,11,12}
{
\draw(v\i)[fill=white] circle(\vr);
}
\end{tikzpicture}
    \caption{Example of a graph constructed by two amalgamations of graphs from $\mathcal{U}$.}
    \label{fig:examples}
\end{figure}
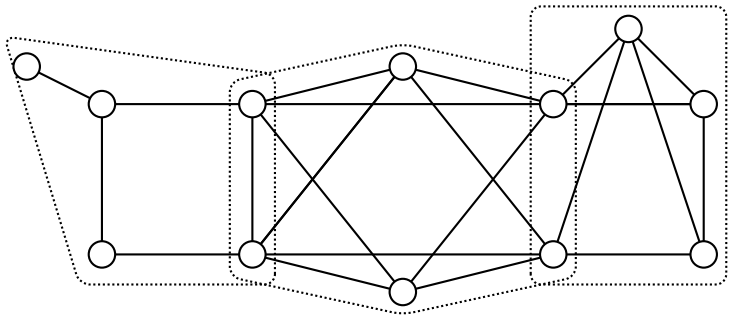

The real power of Theorem~\ref{thm:main} is in applications where the graphs under consideration can be constructed by a repetitive amalgamation of graphs from the family $\mathcal{U}$. Oz and Cangul~\cite{oz-2022} first applied it to  benzenoid chains and in~\cite{oz-2023} extended the method to catacondensed benzenoid systems. In~\cite{chen-2024} the method was elaborated on chains consisting of $8$-cycles. Finally, in~\cite{grad-2025} it was demonstrate that the method can be extended to arbitrary cyclic chains. Theorem~\ref{thm:main} thus widely extends these earlier results. In Fig.~\ref{fig:examples} we give an example which demonstrates that we can get intrinsically different and more complicated examples of this sort. 

%%%%%%%%%%%%%%%%%%%%%%%%%%%%%%
\section{Proof of Theorem~\ref{thm:main}}
\label{sec:proof}

In order to write our transfer matrices compactly, we will use the following abbreviations. First, set 
$$P(G, X, i) = \sum_{x \in X} p(G- x,i)\,.$$ 
Let now $k \ge 1$ be fixed. Then $M(G)$, $M(G,X)$, $0M(G,X)$, and $00M(G,X)$ denote the upper triangular Toeplitz $(k+1) \times (k+1)$ matrices whose first rows are, respectively:  
\begin{align*}
    &  [p(G,0), p(G,1), p(G,2), \dots, p(G,k)]\,, \\
    &  [P(G,X,0), P(G,X,1), P(G,X,2), \dots, P(G,X,k)]\,, \\
    &  [0, P(G,X,0), P(G,X,1), \dots, P(G,X,k-1)]\,, \\
    &  [0, 0, P(G,X,0), P(G,X,1), \dots, P(G,X,k-2)]\,.
\end{align*}

\begin{figure}[ht!]
    \centering
\begin{tikzpicture}[scale=0.8,style=thick,x=1cm,y=1cm]
\def\vr{5pt}

% vertices defined
\coordinate(v1) at (0.0,1.5);
\coordinate(v2) at (0.0,-1.5);

\coordinate(x1) at (-2,2.2);
\coordinate(x2) at (-2,1.95);
\coordinate(x3) at (-2,1.05);

\coordinate(z1) at (-2,0.6);
\coordinate(z2) at (-2,0.35);
\coordinate(z3) at (-2,-0.6);

\coordinate(y1) at (-2,-1.05);
\coordinate(y2) at (-2,-1.3);
\coordinate(y3) at (-2,-2.2);

% Pikice
\node at ($(x2)!0.35!(x3)$) {$\vdots$};
\node at ($(z2)!0.35!(z3)$) {$\vdots$};
\node at ($(y2)!0.35!(y3)$) {$\vdots$};

% \edges                             
\draw (v1) -- (v2); 
\foreach \i in {1,2,3}
{
\draw (x\i) -- (v1);
}
\foreach \i in {1,2,3}
{
\draw (z\i) -- (v1);
\draw (z\i) -- (v2);
}
\foreach \i in {1,2,3}
{
\draw (y\i) -- (v2);
}

% Nariši elipso s poljemi 
\draw (-2,0) ellipse [x radius=0.5cm, y radius=2.5cm];
\draw (-2.5,0.85) -- (-1.5,0.85);
\draw (-2.5,-0.85) -- (-1.5,-0.85);

% Definirajte blob kot krožnico z naključnimi premiki 
\draw[dashed] plot[smooth cycle, tension=.7] coordinates {(-1.3,-2.5) (-1.3,0) 
(-1.3, 2.5) (-3,2.3) (-4,1.5) (-6,2.5) (-6,-2) (-3,-2.4) };

% Risanje polkroga med C in D z uporabo semicircle
\draw (v2) arc[start angle=-120, end angle=120, radius=1.72] -- cycle;
  
% Dodajanje črke S znotraj polkroga
\node at (0:1.1) {\Huge $G'$};
\node at (0:-4.6) {\Huge $H''$};

\node at (-2.8,1.6) {\Large $X$};
\node at (-2.9,0) {\Large $Z$};
\node at (-2.8,-1.6) {\Large $Y$};

%  vertices
\draw[fill=white] (v1) circle (\vr) node[label=above:$x$] {};
\draw[fill=white] (v2) circle (\vr) node[label=below:$y$] {};

\foreach \i in {1,2,3}
{
\draw(x\i)[fill=white, thin] circle(2pt);
}
\foreach \i in {1,2,3}
{
\draw(z\i)[fill=white, thin] circle(2pt);
}
\foreach \i in {1,2,3}
{
\draw(y\i)[fill=white, thin] circle(2pt);
}
\end{tikzpicture}
    \caption{The graph $G$ and notation used in the proof.}
    \label{fig:izrek_ab}
\end{figure}
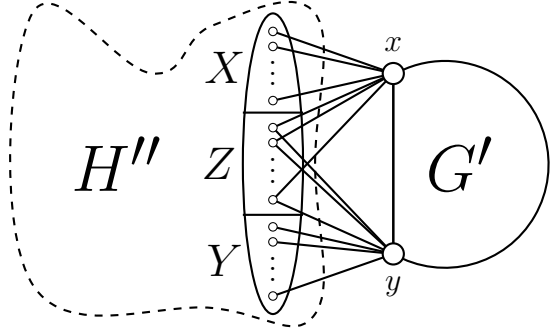

Next let $G\in {\cal U}$ and let $G'$, $x$, $y$, $a$, and $b$ be as in Definition~\ref{def:key}. Then $G-\{x,y\}$ consists of two components $H'$ and $H''$, where $H'$ is a subgraph of $G'$. Define $H$ to be the subgraph of $G$ induced by $V(H'')\cup \{x,y\}$, and let $a,b\in V(H)$, where $\{a,b\}\ne \{x,y\}$, so that one of $a$ and $b$ can be equal to $x$ or $y$. Let further $X = \{x_1, x_2, \dots\}$, $Y = \{y_1, y_2, \dots\}$, and $Z = \{z_1, z_2, \dots\}$ denote the sets of vertices which are adjacent to $x$ but not $y$, to $y$ but not $x$, and to both $x$ and $y$, respectively. See Fig.~\ref{fig:izrek_ab}. 

The proof of Theorem~\ref{thm:main} will consider three cases depending on the  position of vertices $a$ and $b$ related to the neighborhoods of $x$ and $y$. For this sake we also define the set of edges  
$$\mathcal{E} = E_X\cup E_Y\cup E_Z\,,$$ 
where $E_X = \{xu:\ u\in X\}$, $E_Y = \{yu:\ u\in Y\}$, and $E_Z = \{xu:\ u\in Z\} \cup \{yu:\ u\in Z\}$. The first case considered is the following. 

\medskip\noindent 
{\bf Case 1}: $a,b \in V(H'') \setminus (N_G(x) \cup N_G(y))$. \\
Consider the first component of $p_{ab}(G,k)$, that is, $p(G,k)$. We split the $k$-matchings of $G$ into those that contain no edge of $\mathcal{E}$, and into those that contain at least one edge of $\mathcal{E}$. In the latter situation we distinguish the cases when an edge of the $k$-matching lies in $E_X$, in $E_Y$, or in $E_Z$ (where such edges are of two types), and the cases when the $k$-matching contains two edges from $\mathcal{E}$ (there are three typical situations when this can happen).  Setting \mbox{$X = p(G,k)$} and applying Lemmas~\ref{lemma1} and~\ref{lemma2} we can compute as follows:
\begin{align*}
    X = & \ p(G - \bigcup_{e \in \mathcal{E}}e, k)     
    \\   & \ + \sum_{x_i \in X} p(G - x - x_i - \bigcup_{e \in \mathcal{E}} e, k-1) + \sum_{y_j \in Y} p(G - y - y_j - \bigcup_{e \in \mathcal{E}} e, k-1) 
    \\   & \ + \sum_{z_\lambda \in Z} p(G - x - z_\lambda - \bigcup_{e \in \mathcal{E}} e, k-1) + \sum_{z_\lambda \in Z} p(G - y - z_\lambda - \bigcup_{e \in \mathcal{E}} e, k-1) 
      \\   & \ + \sum_{(x_i, y_j) \in (X,Y)} p(G - x - x_i - y - y_j - \bigcup_{e \in \mathcal{E}} e, k-2) 
      \\   & \ + \sum_{(x_i, z_\lambda) \in (X,Z)} p(G - x - x_i - y - z_\lambda - \bigcup_{e \in \mathcal{E}} e, k-2)
      \\   & \ + \sum_{(z_\lambda, y_j) \in (Z,Y)} p(G - x - z_\lambda - y - y_j - \bigcup_{e \in \mathcal{E}} e, k-2)
      \\   = & \ p(G' \cup H'', k) 
      \\  & \ + \sum_{x_i \in X}p((G'-x) \cup (H'' - x_i),k-1) + \sum_{y_j \in Y}p((G'-y) \cup (H'' - y_j),k-1)
      \\  & \ + \sum_{z_\lambda \in Z}(p((G'-x) \cup (H'' - z_\lambda),k-1) + p((G'-y) \cup (H'' - z_\lambda),k-1))
      \\  & \ + \sum_{(x_i, y_j) \in (X,Y)}p((G'-x-y) \cup (H'' - x_i - y_j),k-2)
      \\  & \ + \sum_{(x_i, z_\lambda) \in (X,Z)}p((G'-x-y) \cup (H'' - x_i - z_\lambda),k-2)
      \\  & \ + \sum_{(z_\lambda, y_j) \in (Z,Y)}p((G'-x-y) \cup (H'' - z_\lambda - y_j),k-2)\,.
  \end{align*}
Setting 
\begin{align*}
P(H'', X, k) = & \sum_{x_i \in X} p(H''- x_i,k)\,, 
\end{align*}
and
\begin{align*}
P(H'', (X, Y, Z), k) = & \sum_{(x_i, y_j) \in (X,Y)} p(H''- x_i - y_j,k)\ + \\
& \sum_{(x_i, z_\lambda) \in (X,Z)}p(H'' - x_i - z_\lambda,k)\  + \\
& \sum_{(z_\lambda, y_j) \in (Z,Y)}p(H'' - z_\lambda - y_j,k)\,,
\end{align*}
we can conclude that
\begin{align*}
    X = & \ \Big(p(H'',0), p(H'',1), p(H'',2), \dots ,p(H'',k), 
      \\ & \ 0, P(H'', X, 0) + P(H'', Z, 0), \dots, P(H'', X, k-1) + P(H'', Z, k-1), 
      \\ & \ 0, P(H'', Y, 0) + P(H'', Z, 0), \dots, P(H'', Y, k-1) + P(H'', Z, k-1),  
      \\ & \ 0, 0, P(H'', (X, Y, Z), 0), \dots, P(H'', (X, Y, Z), k-2)\Big) 
      \cdot p_{xy}(G',k)\,.
\end{align*}
The obtained vector corresponds to the first row of the matrix $A = M(H,ab,xy)$ as presented in Fig.~\ref{fig:matrika_ab_zunaj}.

\begin{figure} [ht!]
    \centering
    $$
    \renewcommand{\arraystretch}{2}
    \left(
    \scriptsize
        \begin{array}{c:c:c:c}
            M(H'') & 0M(H'',X \cup Z) & 0M(H'',Y \cup Z) & 00M(H'',(X,Y,Z)) \\
                \hdashline
            M(H''-a) & 0M(H''-a,X \cup Z) & 0M(H''-a,Y \cup Z) & 00M(H''-a,(X,Y,Z)) \\
                \hdashline
            M(H''-b) & 0M(H''-b,X \cup Z) & 0M(H''-b,Y \cup Z) & 00M(H''-b,(X,Y,Z)) \\
                \hdashline
            M(H''-a-b) & 0M(H''-a-b,X \cup Z) & 0M(H''-a-b,Y \cup Z) & 00M(H''-a-b,(X,Y,Z)) \\
        \end{array}
    \right)
    $$
    \caption{Matrix $A = M(G,ab,xy)$ of Case 1.}
    \label{fig:matrika_ab_zunaj}
\end{figure}

The above proof similarly holds for the remaining $k$ components of the form $p(G, i)$, where $0 \leq i \leq k-1$. In this way we obtain the four upper triangular Toeplitz matrices of dimension $(k + 1) \times (k + 1)$ which lie above in the matrix $M(H,ab,xy)$.

Since $a, b \notin N_G(x) \cup N_G(y)$, the computation of the components $p(G-a,k)$, $p(G-b,k)$, and $p(G-a-b,k)$ proceeds along the same lines. In this way we arrive at an analogous result, the only difference being that the graph $H''$ is accordingly altered to $H'' - a$, $H'' - b$, and $H''-a-b$. Thus, we obtain the remaining parts of the matrix $A = M(H,ab,xy)$ from Fig.~\ref{fig:matrika_ab_zunaj}.

\medskip\noindent
{\bf Case 2}: $a,b \in V(H'')$, $|\{x, y, a, b\}| = 4$. \\
The subcase when $a,b \notin N_G(x) \cup N_G(y)$ is covered by Case 1, hence we can assume that at least one of $a$ and $b$ lies in $N_G(x)\cup N_G(y)$. Without loss of generality assume that $a \in X$ and that $a = x_1$. The first component $p(G,k)$ remains the same as in Case 1 because the vertex $a$ does not influence it. This observation extends to the first $k+1$ components and even to the components $p(G-b,i)$ for some $i \in [k]$ later on.

The component $p(G-a,k)$ changes only in the summation over the set $X$ which is replaced by $X \setminus \{x_1\}$ or equivalently $X \setminus \{a\}$. The same holds for the summation over sets $(X,Y)$ and $(X,Z)$. The component $p(G-a-b,k)$ changes in the same manner as $p(G-a,k)$.

Therefore, by identifying the position of vertices $a$ and $b$ relative to sets $X$, $Y$ and $Z$, we can construct a matrix $M(H,ab,xy)$, such as the one in Fig.~\ref{fig:matrika_ab_zunaj}, with minor modifications. 

\medskip\noindent
{\bf Case 3}: $a,b \in V(H'')$, $\{a,b\} \neq \{x,y\}$. \\
Due to Case 2 it suffices to consider the case $|\{x, y, a, b\}| = 3$. Assume first that $a = x$. The first component $p(G,k)$ remains the same as in Case 1, as the vertex $a$ does not influence it. This observation extends to the first $k+1$ components and even to components $p(G-b,i)$ for some $i \in [k]$ later on. We now proceed to compute $p(G-a, k)$. Setting $X_a = p(G-a,k)$ and $P(H'', Y, k) = \sum_{y_j \in Y} p(H''- y_j,k)$ we can compute as follows:

  \begin{align*}
    X_a = & \ p(G - a - \bigcup_{e \in \mathcal{E}}e, k)     
    \\   & \  + \sum_{y_j \in Y} p(G - a - y - y_j - \bigcup_{e \in \mathcal{E}} e, k-1) 
    \\   & \  + \sum_{z_\lambda \in Z} p(G - a - y - z_\lambda - \bigcup_{e \in \mathcal{E}} e, k-1) 
      \\   = & \ p((G' - x) \cup H'', k) 
      \\  & \  + \sum_{y_j \in Y}p((G' - x -y) \cup (H'' - y_j),k-1)
      \\  & \ + \sum_{z_\lambda \in Z}p((G'-x-y) \cup (H''-z_\lambda),k-1)
      \\ = & \ \Big( 0,\quad  0,\quad  0, \quad \dots, \quad 0,
      \\ & \ p(H'',0), p(H'',1), p(H'',2), \dots ,p(H'',k), 
       \\ & \ 0,\quad  0,\quad  0, \quad \dots, \quad 0,
      \\ & \ 0, P(H'', Y, 0) + P(H'', Z, 0), \dots, P(H'', Y, k-1) + P(H'', Z, k-1)\Big) 
      \cdot p_{xy}(G',k)\,.
  \end{align*}

The computation is similar for all components $p(G-a,i)$ and $p(G-a-b, i)$ for some $i \in [k]$. In this way we have computed the matrix $A = M(H,ab,xy)$ for the case where $a = x$, it is presented in Fig.~\ref{fig:matrika_ax}.  

\begin{figure} [ht!]
    \centering
    $$
    \renewcommand{\arraystretch}{1.9}
    \left(
    \scriptsize
        \begin{array}{c:c:c:c}
            M(H'') & 0M(H'',X \cup Z) & 0M(H'',Y \cup Z) & 00M(H'',(X,Y,Z)) \\
                \hdashline
            0 & M(H'') & 0 & 0M(H'',Y \cup Z) \\
                \hdashline
            M(H''-b) & 0M(H''-b,X \cup Z) & 0M(H''-b,Y \cup Z) & 00M(H''-b,(X,Y,Z)) \\
                \hdashline
            0 & M(H''-b) & 0 & 0M(H''-b,Y \cup Z)
        \end{array}
    \right)
    $$
    \caption{The matrix $A = M(G,ab,xy)$ for the case $a = x$.}
    \label{fig:matrika_ax}
\end{figure}

The cases when $a=y$ or $b \in \{x,y\}$ are treated similarly, the corresponding matrices are respectively presented in Figs.~\ref{fig:matrika_ay}, \ref{fig:matrika_bx}, and \ref{fig:matrika_by}. Considering all the cases  on the position of the vertices $a$ and $b$, the proof of Theorem~\ref{thm:main} is completed. 

\begin{figure} [ht!]
    \centering
    $$
    \renewcommand{\arraystretch}{1.9}
    \left(
    \scriptsize
        \begin{array}{c:c:c:c}
            M(H'') & 0M(H'',X \cup Z) & 0M(H'',Y \cup Z) & 00M(H'',(X,Y,Z)) \\
                \hdashline
            0 & 0 & M(H'') & 0M(H'',X \cup Z) \\
                \hdashline
            M(H''-b) & 0M(H''-b,X \cup Z) & 0M(H''-b,Y \cup Z) & 00M(H''-b,(X,Y,Z)) \\
                \hdashline
            0 & 0 & M(H''-b) & 0M(H''-b,X \cup Z)
        \end{array}
    \right)
    $$
    \caption{The matrix $A = M(G,ab,xy)$ for the case $a = y$.}
    \label{fig:matrika_ay}
\end{figure}

\begin{figure} [ht!]
    \centering
    $$
    \renewcommand{\arraystretch}{1.9}
    \left(
    \scriptsize
        \begin{array}{c:c:c:c}
            M(H'') & 0M(H'',X \cup Z) & 0M(H'',Y \cup Z) & 00M(H'',(X,Y,Z)) \\
                \hdashline
            M(H''-a) & 0M(H''-a,X \cup Z) & 0M(H''-a,Y \cup Z) & 00M(H''-a,(X,Y,Z)) \\
                \hdashline
            0 & M(H'') & 0 & 0M(H'',Y \cup Z) \\
                \hdashline
            0 & M(H''-a) & 0 & 0M(H''-a,Y \cup Z)
        \end{array}
    \right)
    $$
    \caption{The matrix $A = M(G,ab,xy)$ for the case $b = x$.}
    \label{fig:matrika_bx}
\end{figure}

\begin{figure} [ht!]
    \centering
    $$
    \renewcommand{\arraystretch}{1.9}
    \left(
    \scriptsize
        \begin{array}{c:c:c:c}
            M(H'') & 0M(H'',X \cup Z) & 0M(H'',Y \cup Z) & 00M(H'',(X,Y,Z)) \\
                \hdashline
            M(H''-a) & 0M(H''-a,X \cup Z) & 0M(H''-a,Y \cup Z) & 00M(H''-a,(X,Y,Z)) \\
                \hdashline
            0 & 0 & M(H'') & 0M(H'',X \cup Z) \\
                \hdashline
            0 & 0 & M(H''-a) & 0M(H''-a,X \cup Z)
        \end{array}
    \right)
    $$
    \caption{The matrix $A = M(G,ab,xy)$ for the case $b = y$.}
    \label{fig:matrika_by}
\end{figure}

%%%%%%%%%%%%%%%%%%%%%%%%%%%%%%
\section{Two examples}
\label{sec:2-examples}
%%%%%%%%%%%%%%%%%%%%%%%%%%%%%%

To demonstrate how the method developed in this paper works, we present in this section two examples.

\subsection{A sporadic example}

\begin{figure}[ht!]
  \centering
  \begin{tikzpicture}[scale=0.8,style=thick,x=1cm,y=1cm]
    \def\vr{5pt}
    % Prvi trikotnik
    \node[draw, circle, inner sep=2.5pt, fill=black, label=below:$a$] (v1) at (-2.732/2,-1/2) {};
    \node[draw, circle, inner sep=2.5pt, fill=black, label=above:$b$] (v2) at (-2.732/2,1/2) {};
    \node[draw, circle, inner sep=2.5pt, fill=black] (v3) at (-1/2,0/2) {};
    \draw (v1) -- (v2) -- (v3) -- (v1);
    
    % Drugi trikotnik
    \node[draw, circle, inner sep=2.5pt, fill=black] (v6) at (1/2,0) {};
    \node[draw, circle, inner sep=2.5pt, fill=black, label=below:$x$] (v4) at (2.732/2,-1/2) {};
    \node[draw, circle, inner sep=2.5pt, fill=black, label=above:$y$] (v5) at (2.732/2,1/2) {};
    \draw (v6) -- (v4) -- (v5) -- (v6);
    
    % Povezava med vozliščema
    \draw (v3) -- (v6);

    % Šesterokotnik pritrjen na AB
    \node[draw, circle, inner sep=2.5pt, fill=black] (v10) at (-4.464/2, -2/2) {};
    \node[draw, circle, inner sep=2.5pt, fill=black, label = left:$v_1$] (v7) at (-6.19615/2, -1/2) {};
    \node[draw, circle, inner sep=2.5pt, fill=black, label = left:$v_2$] (v8) at (-6.19615/2, 1/2) {};
    \node[draw, circle, inner sep=2.5pt, fill=black] (v9) at (-4.464/2, 2/2) {};
    \draw (v1) -- (v10) -- (v7) -- (v8) -- (v9) -- (v2);

    % Petkotnik pritrjen na XY
    \node[draw, circle, inner sep=2.5pt, fill=black, label =below right:$v_3$] (v12) at (4.464/2, -1.7/2) {};
    \node[draw, circle, inner sep=2.5pt, fill=black, label =right:$v_4$] (v11) at (5.8/2, 0) {};
    \node[draw, circle, inner sep=2.5pt, fill=black] (v13) at (4.464/2, 1.7/2) {};
    \draw (v4) -- (v12) -- (v11) -- (v13) -- (v5);
% Vozlišča
    \foreach \i in {1,2,3,4,5,6,7,8,9,10,11,12,13}
{
\draw(v\i)[fill=white] circle(\vr);
}
\end{tikzpicture}
  \caption{The graph $G$.}
  \label{fig:trikotmol}
\end{figure}
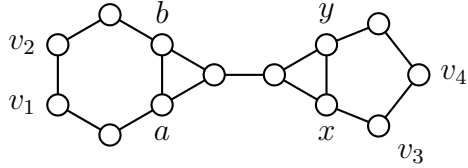

In the  first example we compute the Hosoya index of the graph $G$ from Fig.~\ref{fig:trikotmol}.

\begin{figure}[h!]
  \centering
  \begin{tikzpicture}[scale=1.0,style=thick,x=1cm,y=1cm]
    \def\vr{5pt}
    % Prvi trikotnik
    \node[draw, circle, inner sep=2.5pt, fill=black, label=left:$a$] (v1) at (-2.732,-1) {};
    \node[draw, circle, inner sep=2.5pt, fill=black, label=left:$b$] (v2) at (-2.732,1) [circle,draw] {};
    \node[draw, circle, inner sep=2.5pt, fill=black] (v3) at (-1,0) [circle,draw] {};
    \draw (v1) -- (v2) -- (v3) -- (v1);
    
    % Drugi trikotnik
    \node[draw, circle, inner sep=2.5pt, fill=black, label=above left:$z$] (v4) at (1,0) [circle,draw] {};
    \node[draw, circle, inner sep=2.5pt, fill=black, label=right:$x$] (v5) at (2.732,-1) [circle,draw] {};
    \node[draw, circle, inner sep=2.5pt, fill=black, label=right:$y$] (v6) at (2.732,1) [circle,draw] {};
    \draw (v4) -- (v5) -- (v6) -- (v4);
    
    % Povezava med vozliščema
    \draw (v3) -- (v4);

    \foreach \i in {1,2,3,4,5,6}
{
\draw(v\i)[fill=white] circle(\vr);
}
\end{tikzpicture}
  \caption{The graph $T$.}
  \label{fig:trikot}
\end{figure}
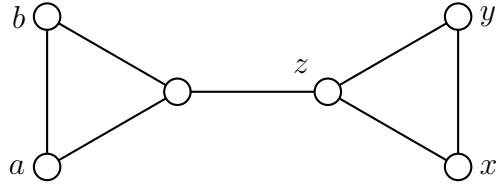

For our purposes, the graph $G$ can be decomposed into a hexagon (left), a graph $T$ as shown in Fig.~\ref{fig:trikot}, and a pentagon (right). 

Consider first the graph $T$, and let $x$, $y$, $a$, and $b$ be its vertices as marked in Fig.~\ref{fig:trikot}. According to the setup of the proof of Theorem~\ref{thm:main}, Case 1 of the proof applies to this situation. Note that using the notation from the proof, $X = Y = \emptyset$ and $Z = \{z\}$. Setting $T'' = T - x -y$, the transfer matrix $M(T,ab,xy)$ can be written as follows:

$$
\renewcommand{\arraystretch}{2}
\left(
\scriptsize
    \begin{array}{c:c:c:c}
        M(T'') & 0M(T'', \{z\}) & 0M(T'',\{z\}) & 00M(T'',(\emptyset,\emptyset,\{z\})) \\
            \hdashline
        M(T''-a) & 0M(T''-a, \{z\}) & 0M(T''-a, \{z\}) & 00M(T''-a,(\emptyset,\emptyset,\{z\})) \\
            \hdashline
        M(T''-b) & 0M(T''-b, \{z\}) & 0M(T''-b, \{z\}) & 00M(T''-b,(\emptyset,\emptyset,\{z\})) \\
            \hdashline
        M(T''-a-b) & 0M(T''-a-b, \{z\}) & 0M(T''-a-b, \{z\}) & 00M(T''-a-b,(\emptyset,\emptyset,\{z\})) \\
    \end{array}
\right)
$$
Clearly, $T'' - a \cong T'' - b \cong  P_3$ and  $T'' - a - b \cong P_2$. Note further that $M(G,(\emptyset,\emptyset,Z)) = 0$ and $M(G,\{z\})$ = $M(G-z)$. Using this facts we can write:

% Če bi želeli vmesno matriko pri izračunu.
% $$
% M(T,ab,xy) =
% \renewcommand{\arraystretch}{2}
% \left(
% \scriptsize
%     \begin{array}{c:c:c:c}
%         M(T') & 0M(T', \{z\}) & 0M(T',\{z\}) & 0 \\
%             \hdashline
%         M(P_3) & 0M(P_3, \{z\}) & 0M(P_3, \{z\}) & 0 \\
%             \hdashline
%         M(P_3) & 0M(P_3, \{z\}) & 0M(P_3, \{z\}) & 0 \\
%             \hdashline
%         M(P_3) & 0M(P_2, \{z\}) & 0M(P_2, \{z\}) & 0 \\
%     \end{array}
% \right)
% $$

$$
M(T,ab,xy) =
\renewcommand{\arraystretch}{2}
\left(
\scriptsize
    \begin{array}{c:c:c:c}
        M(T'') & 0M(K_3) & 0M(K_3) & 0 \\
            \hdashline
        M(P_3) & 0M(P_2) & 0M(P_2) & 0 \\
            \hdashline
        M(P_3) & 0M(P_2) & 0M(P_2) & 0 \\
            \hdashline
        M(P_2) & 0M(P_1) & 0M(P_1) & 0 \\
    \end{array}
\right)
$$
% Given that the largest graph which appears within the adjacency matrix $M(T,ab,xy)$ is $T''$, it suffices to compute just $1$ and $2$-matchings. 
Setting the first row of the matrices $M(a,b,c)$, $0M(a,b,c)$, and $00M(a,b,c)$ as 
\begin{align*}
    & [a, b, c, 0, \dots, 0]\,, \\
    & [0, a, b, c, 0, \dots, 0]\,, \\
    & [0, 0, a, b, c, 0, \dots, 0]\,,
\end{align*}
the matrix $M(T,ab,xy)$ can be written as follows: 
$$
M(T,ab,xy) = 
\renewcommand{\arraystretch}{2}
    \left(
        \scriptsize
        \begin{array}{c:c:c:c}
            M(1,4,1) & 0M(1,3,0) & 0M(1,3,0) & 0 \\
                \hdashline
            M(1,2,0) & 0M(1,1,0) & 0M(1,1,0) & 0 \\
                \hdashline
            M(1,2,0) & 0M(1,1,0) & 0M(1,1,0) & 0\\
                \hdashline
            M(1,1,0) & 0M(1,0,0) & 0M(1,0,0) & 0 \\
        \end{array}
    \right)\,.
$$

Following a similar approach, we construct the transfer matrices for the hexagon and pentagon.

$$
M(C_6,v_1v_2,ab) = 
\renewcommand{\arraystretch}{2}
    \left(
        \scriptsize
        \begin{array}{c:c:c:c}
            M(1,3,1) & 0M(1,2,0) & 0M(1,2,0) & 00M(1,1,0) \\
                \hdashline
            M(1,1,0) & 0M(1,1,0) & 0M(1,0,0) & 00M(1,0,0) \\
                \hdashline
            M(1,1,0) & 0M(1,0,0) & 0M(1,1,0) & 00M(1,0,0) \\
                \hdashline
            M(1,0,0) & 0M(1,0,0) & 0M(1,0,0) & 00M(1,0,0) \\
        \end{array}
    \right)\,,
$$

$$
M(C_5,xy,v_3v_4) = 
\renewcommand{\arraystretch}{2}
    \left(
        \scriptsize
        \begin{array}{c:c:c:c}
            M(1,2,0) & 0M(1,1,0) & 0M(1,1,0) & 00M(1,0,0) \\
                \hdashline
            M(1,1,0) & 0 & 0M(1,0,0) & 0 \\
                \hdashline
            M(1,0,0) &  0M(1,0,0) &  0M(1,0,0) & 00M(1,0,0)\\
                \hdashline
            M(1,0,0) & 0 & 0M(1,0,0) & 0 \\
        \end{array}
    \right)\,.
$$
We note that the latter two matrices can also be deduced from~\cite[[Theorem 1]{grad-2025}. Since we have decomposed $G$ into a hexagon, $T$, and pentagon, we obtain:
\begin{align*} 
p_{v_1v_2}(G,6) &= M(C_6,v_1v_2,ab) \cdot M(T,ab,xy) \cdot M(C_5,xy,v_3v_4) \cdot p_{v_3v_4}(P_2,6) \\
    &= [26,\ 177,\ 336,\ 261,\ 95,\ 16,\ 1,\ \dots] \,,
\end{align*}
and conclude from here that 
$$Z(G) = 26 + 177 + 336 + 261 + 95 + 16 + 1 = 912\,.$$

%%%%%%%%%%%%%%%%%%%%%%%%%%%%%%
\subsection{A chemical example}
\label{sec:chem-example}
%%%%%%%%%%%%%%%%%%%%%%%%%%%%%%

As the second example of how the method developed in this paper applies to chemical graphs, consider the molecular graph $G$ from Fig.~\ref{fig:molecule}. In this way we demonstrate that the developed method allows us to deal also with certain classes of pericondensed benzenoid systems as well as other non cyclic structures.

\begin{figure}[ht!]
    \centering
    \begin{tikzpicture}[scale=0.8,style=thick,x=1cm,y=1cm]

\def\vr{5pt}
% vertices defined
\coordinate(v1) at (-5,0.5);
\coordinate(v2) at (-5,-0.5);
\coordinate(v3) at (-4.0,2.0);
\coordinate(v4) at (-4.0,1.0);
\coordinate(v5) at (-4,-1.0);
\coordinate(v6) at (-3,2.5);
\coordinate(v7) at (-3.0,0.5);
\coordinate(v8) at (-3.0,-0.5);
\coordinate(v9) at (-2.0,2);
\coordinate(v10) at (-2.0,1.0);
\coordinate(v11) at (-2.0,-1.0);
\coordinate(v12) at (-1.0,0.5);
\coordinate(v13) at (-1.0,-0.5);
\coordinate(v14) at (0,0.75);
\coordinate(v15) at (0,-0.75);
\coordinate(v16) at (0.75,0);
\coordinate(v17) at (2.25,0);
\coordinate(v18) at (2.5,1.5);
\coordinate(v19) at (3.0,0.75);
\coordinate(v20) at (3.0,-0.75);
\coordinate(v21) at (4,0.5);
\coordinate(v22) at (4,-0.5);
% \edges                             
\draw (v1) -- (v2) -- (v5) -- (v8) -- (v7) -- (v4) -- cycle; 
\draw (v4) -- (v3) -- (v6) -- (v9) -- (v10) -- (v7) -- cycle; 
\draw (v7) -- (v10) -- (v12) -- (v13) -- (v11) -- (v8) -- cycle; 
\draw (v12) -- (v14) -- (v16) -- (v15) -- (v13) -- cycle; 
\draw (v16) -- (v17);
\draw (v18) -- (v19);
\draw (v17) -- (v19) -- (v21) -- (v22) -- (v20) -- cycle; 

%  vertices
\foreach \i in {1,2,3,4,5,6,7,8,9,10,11,12,13,14,15,16,17,18,19,20,21,22}
{
\draw(v\i)[fill=white] circle(\vr);
}
\end{tikzpicture}
    \caption{Graph $G$.}
    \label{fig:molecule}
\end{figure}

\begin{figure}[ht!]
    \centering
    \begin{tikzpicture}[scale=0.8,style=thick,x=1cm,y=1cm]
\def\vr{5pt}
% vertices defined
\coordinate(v1) at (-5,0.5);
\coordinate(v2) at (-5,-0.5);
\coordinate(v3) at (-4.0,2.0);
\coordinate(v4) at (-4.0,1.0);
\coordinate(v5) at (-4,-1.0);
\coordinate(v6) at (-3,2.5);
\coordinate(v7) at (-3.0,0.5);
\coordinate(v8) at (-3.0,-0.5);
\coordinate(v9) at (-2.0,2);
\coordinate(v10) at (-2.0,1.0);
\coordinate(v11) at (-2.0,-1.0);
\coordinate(v12) at (-1.0,0.5);
\coordinate(v13) at (-1.0,-0.5);
\coordinate(v14) at (0,0.75);
\coordinate(v15) at (0,-0.75);
\coordinate(v16) at (0.75,0);
\coordinate(v17) at (2.25,0);
\coordinate(v18) at (2.5,1.5);
\coordinate(v19) at (3.0,0.75);
\coordinate(v20) at (3.0,-0.75);
\coordinate(v21) at (4,0.5);
\coordinate(v22) at (4,-0.5);
% \edges                             
\draw (v1) -- (v2) -- (v5) -- (v8) -- (v7) -- (v4) -- cycle; 
\draw (v4) -- (v3) -- (v6) -- (v9) -- (v10) -- (v7) -- cycle; 
\draw (v7) -- (v10) -- (v12) -- (v13) -- (v11) -- (v8) -- cycle; 
\draw (v12) -- (v14) -- (v16) -- (v15) -- (v13) -- cycle; 
\draw (v16) -- (v17);
\draw (v18) -- (v19);
\draw (v17) -- (v19) -- (v21) -- (v22) -- (v20) -- cycle; 

\draw[rounded corners, densely dotted](-3,3)--(-5.3,1.6)--(-5.3,-1.3)--(-0.7,-1.3)--(-0.7,1)--(-0.7,1.6)--cycle;
\draw[rounded corners, densely dotted](-1.3,-1)--(-1.3,1)--(1,1)--(2.75,0)--(1,-1)--cycle;
\draw[rounded corners, densely dotted](1.7,-1)--(0.25,0)--(1.7,1)--(2.3,2)--(4.3,2)--(4.3,-1)--cycle;
%  vertices
\foreach \i in {1,2,3,4,5,6,7,8,9,10,11,12,13,14,15,16,17,18,19,20,21,22}
{
\draw(v\i)[fill=white] circle(\vr);
}
% Add labels to specific vertices
\node[left=5pt] at (v1) {$a$};
\node[left= 5pt] at (v2) {$b$};
\node[above left = 1pt] at (v4) {$c$}; 
\node[below = 2pt] at (v8) {$d$};
\node[above right = 1pt] at (v10) {$e$};
\node[above=10pt] at (v12) {$f$};
\node[below = 1pt] at (v13) {$g$};
\node[below =7pt] at (v16) {$h$};
\node[below =7pt] at (v17) {$i$};
\node[right =5pt] at (v21) {$x$};
\node[right =5pt] at (v22) {$y$};
\end{tikzpicture}
    \caption{Decomposition of graph $G$.}
    \label{fig:de1}
\end{figure}
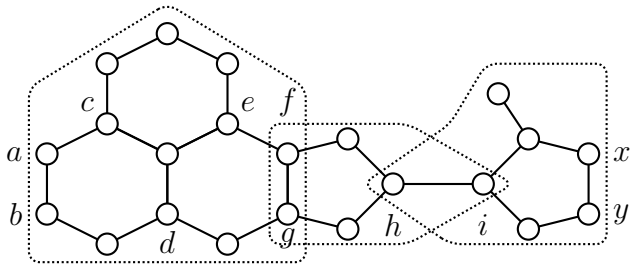

The graph $G$ can be decomposed in many different ways such that Theorem~\ref{thm:main} can be applied. One possible decomposition is shown in Fig.~\ref{fig:de1}.

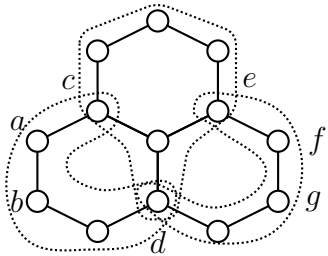
\begin{figure}[ht!]
    \centering
    \begin{tikzpicture}[scale=0.8,style=thick,x=1cm,y=1cm]
\def\vr{5pt}
% vertices defined
\coordinate(v1) at (-5,0.5);
\coordinate(v2) at (-5,-0.5);
\coordinate(v3) at (-4.0,2.0);
\coordinate(v4) at (-4.0,1.0);
\coordinate(v5) at (-4,-1.0);
\coordinate(v6) at (-3,2.5);
\coordinate(v7) at (-3.0,0.5);
\coordinate(v8) at (-3.0,-0.5);
\coordinate(v9) at (-2.0,2);
\coordinate(v10) at (-2.0,1.0);
\coordinate(v11) at (-2.0,-1.0);
\coordinate(v12) at (-1.0,0.5);
\coordinate(v13) at (-1.0,-0.5);

% \edges                        
\draw (v1) -- (v2) -- (v5) -- (v8) -- (v7) -- (v4) -- cycle; 
\draw (v4) -- (v3) -- (v6) -- (v9) -- (v10) -- (v7) -- cycle; 
\draw (v7) -- (v10) -- (v12) -- (v13) -- (v11) -- (v8) -- cycle; 
\draw[densely dotted] plot[smooth cycle, tension=.8] coordinates {(-4,1.3)(-5.3,0.6) (-5.3,-1) (-4,-1.3) (-3,-1) (-2.8,-0.2) (-4,-0.3) (-4.5,0.2) (-3.7,0.8)};
\draw[densely dotted] plot[smooth cycle, tension=.8] coordinates {(-0.8,0.8) (-0.8,-0.7) (-2, -1.2) (-3.2, -0.6) (-3.1, -0.2) (-2,-0.6) (-1.2,0) (-2,0.7) (-2.4, 1) (-2,1.3)};
\draw[rounded corners, densely dotted](-4.3,2.3)--(-4.3,1)--(-3.5,0.4)--(-3.3,-0.8)--(-2.7,-0.8)--(-2.5,0.4)--(-1.7,1)--(-1.7,2.3)--(-3,2.8)--cycle;
\foreach \i in {1,2,3,4,5,6,7,8,9,10,11,12,13}
{
\draw(v\i)[fill=white] circle(\vr);
}
%Labels
\node[above left] at (v1) {$a$};
\node[left=1pt] at (v2) {$b$};
\node[above left = 6pt] at (v4) {$c$}; 
\node[below = 7pt] at (v8) {$d$};
\node[above right = 7pt] at (v10) {$e$};
\node[right=7pt] at (v12) {$f$};
\node[right=6pt] at (v13) {$g$};
% \node[below =7pt] at (v16) {$h$};
% \node[below =7pt] at (v17) {$i$};
% \node[right =5pt] at (v21) {$x$};
% \node[right =5pt] at (v22) {$y$};
\end{tikzpicture}
    \caption{Decomposition of phenalene.}
    \label{fig:de2}
\end{figure}

The decomposition in Fig.~\ref{fig:de1} may appear to be the most natural one, but computing the transfer matrix for the phenalene $P$, that is, the left graph from the figure's decomposition, can be challenging. This is why we consider a smaller further decomposition of $P$ as presented in Fig.~\ref{fig:de2}. This indicates that the main advantage and motivation behind our method is the ability to break down the graph into smaller parts.

Given that the graph $G$ contains 22 vertices, we know that its matching number is at most $11$. Therefore, we set $k = 11$. In the decomposition from Fig.~\ref{fig:de2}, we calculate the transfer matrix of the phenalene $P$ as follows:
$$
M(P,ab,fg) = M(P_5, ab, cd) \cdot M(T, cd, ed) \cdot M(P_5, ed, fg)\ ,
$$
where $T$ is the graph of the toluene molecule. Denoting the molecular graph of the methylcyclopentane by $F_1$ and the molecular graph of the 1,2-dimethylcyclopentene by $F_2$, we can compute as follows:
\begin{align*} 
p_{ab}(G,11) &= M(P,ab,fg) \cdot M(F_1,fg,hi) \cdot M(F_2,hi,xy) \cdot  p_{xy}(P_2,11) \\
    &= [3, 182, 2098, 9036, 18676, 21476, 14867, 6425, 1740, 286, 26, 1]\,.
\end{align*}
We can conclude that 
$$Z(G) = \sum_{i=1}^{k+1} (p_{ab}(G,11))_i = 74816\,.$$

%%%%%%%%%%%%%%%%%%%%%%%%%%%%%%%%%%%%%%%%%%%%%%%%%%%%%%%%
\section*{Acknowledgements}
%%%%%%%%%%%%%%%%%%%%%%%%%%%%%%%%%%%%%%%%%%%%%%%%%%%%%%%%

and S.\ Klav\v{z}ar has been supported by the Slovenian Research Agency ARIS (research core funding P1-0297 and projects N1-0285, N1-0355). 

%%%%%%%%%%%%%%%%%%%%%%%%%%%%
\section*{Declaration of interests}
%%%%%%%%%%%%%%%%%%%%%%%%%%%%
 
The authors declare that they have no conflict of interest. 

%%%%%%%%%%%%%%%%%%%%%%%%%%%%
\section*{Data availability}
%%%%%%%%%%%%%%%%%%%%%%%%%%%%
 
Our manuscript has no associated data.

%%%%%%%%%%%%%%%%%%%%%%%%%%%%%%%%%%%%%%%%%%%%%%

\end{document}